\documentstyle[12pt, fullpage]{amsart}
\begin{document}
\title{The Smallest Solution of $\phi(30n+1)<\phi(30n)$ is
\dots}\author{Greg Martin}\maketitle

\def\mod#1{{\ifmmode\text{\rm\ (mod~$#1$)}
 \else\discretionary{}{}{\hbox{ }}\rm(mod~$#1$)\fi}}

In a previous issue of the {\it American Mathematical Monthly\/},
D.~J.~Newman \cite{newman} showed that for any positive integers $a$,
$b$, $c$, and $d$ with $ad\ne bc$, there exist infinitely many
positive integers $n$ for which $\phi(an+b)<\phi(cn+d)$, where
$\phi(m)$ is the familiar Euler totient function, the number of
positive integers less than and relatively prime to $m$. In
particular, it must be the case that $\phi(30n+1)<\phi(30n)$
infinitely often; however, Newman mentions that there are no solutions
of this inequality with $n\le{}$20,000,000, and he states that a
solution ``is not explicitly available and it may be beyond the reach
of any possible computers''. The purpose of this note is to describe a
method for computing solutions to inequalities of this type that
avoids the need to factor large numbers. In particular, we explicitly
compute the smallest number $n$ satisfying $\phi(30n+1)<\phi(30n)$.

It is quite easy to compute values of $n$ for which $\phi(30n+1)$ is
relatively small by imposing many congruence conditions on $n$ modulo
primes, so that $30n+1$ is highly composite. However, the numbers $n$
that arise in this way are quite large, having hundreds of
digits. Computing $\phi(30n)$ exactly relies on the factorization of
$30n$, which for integers of this size is not possible to find in a
reasonable amount of time with today's computers and factoring
algorithms. The idea underlying our method is to use partial knowledge
of the factorization of a large number $m$ to get an estimate for
$\phi(m)$. We rely on the following claim:

\vskip12pt
{\noindent\bf Claim 1. \it Let $p_i$ denote the $i^{th}$ prime
number. Let $q=\prod_{i=r+1}^{r+s}p_i$ for some positive integers $r$
and $s$, and let $m$ be an integer that is not divisible by any of the
primes $p_1$, \dots, $p_r$. Then:

{\rm(a)} if $m\le q$, then $m$ has at most $s$ distinct prime
factors;

{\rm(b)} if $m$ has at most $s$ distinct prime factors, then
$\phi(m)/m \ge \phi(q)/q$.}
\vskip12pt

{\noindent\it Proof.\/} Let $t$ be the number of distinct prime
factors of $m$, and let the prime factors be $p_{i_1}$, \dots,
$p_{i_t}$ with $i_1<\dots<i_t$. Since none of the primes $p_1$, \dots,
$p_r$ divide $m$, it must be the case that $i_1\ge r+1$, $i_2\ge r+2$,
and so on. So if we define $k=\prod_{j=r+1}^{r+t}p_j$, we see that
$k\le\prod_{j=1}^tp_{i_j}\le m$. But $m\le q$ by assumption, and so
$k\le q$, which can clearly only be the case if $t\le s$. This proves
part (a) of the claim.

As for part (b), we use the fact that the function $\phi(m)/m$ can be
written as a product over primes dividing $m$:
\begin{equation*}
{\phi(m)\over m} = \prod_{p\mid m} \big( 1-\frac1p \big).
\end{equation*}
With $k$ defined as above, notice that
\begin{equation*}
{\phi(m)\over m} = \prod_{j=1}^t \big( 1-\frac1{p_{i_j}} \big) \ge
\prod_{j=1}^t \big( 1-\frac1{p_{r+j}} \big) = {\phi(k)\over k},
\end{equation*}
since $1-1/p$ is an increasing function of $p$. On the other hand,
since $t\le s$ by assumption, we have
\begin{equation*}
{\phi(k)\over k} = \prod_{j=r+1}^{r+t} \big( 1-\frac1{p_j} \big) \ge
\prod_{j=r+1}^{r+s} \big( 1-\frac1{p_j} \big) = {\phi(q)\over q},
\end{equation*}
since each $1-1/p$ is less than 1. This proves part (b) of the claim.
\vskip12pt

We now proceed to find the smallest solution of
$\phi(30n+1)<\phi(30n)$, though it must be pointed out that the method
applies to any inequality of the form $\phi(an+b)<\phi(cn+d)$. Clearly
$30n+1\equiv1\mod{30}$ no matter what $n$ is. Also, if $n$ is a
solution of $\phi(30n+1)<\phi(30n)$, then we must have
\begin{equation*}
{\phi(30n+1)\over30n+1} < {\phi(30n)\over30n+1} < {\phi(30)n\over30n}
= \frac4{15} = 0.26666\dots,
\end{equation*}
since the inequality $\phi(ab)\le\phi(a)b$ holds for any numbers $a$
and $b$. Thus it makes sense to look for numbers that satisfy both
these conditions.

\vskip12pt
{\noindent\bf Claim 2. \it Let $z=(p_4p_5\dotsm
p_{383})p_{385}p_{388}$. Then $z$ is the smallest positive integer
satisfying $z\equiv1\mod{30}$ and $\phi(z)/z<4/15$.}
\vskip12pt

{\noindent\it Proof.\/} A computation shows that $z$ is indeed
congruent to $1\mod{30}$ and that
\begin{equation*}
\frac{\phi(z)}z = \bigg( \prod_{i=4}^{383} \big( 1-\frac1{p_i} \big)
\bigg) \big( 1-\frac1{p_{385}} \big) \big( 1-\frac1{p_{388}} \big) =
0.2666117\ldots<\frac4{15}.
\end{equation*}
Suppose $m$ is an integer satisfying $m\equiv1\mod{30}$ and
$\phi(m)/m<4/15$. Because of the congruence condition, $m$ cannot be
divisible by 2, 3, or 5. If we define $q_1=\prod_{i=4}^{384}p_i$, then
we can compute that $\phi(q_1)/q_1 = 0.26671\dots$, and so
$\phi(q_1)/q_1>\phi(m)/m$. Thus if we apply part (b) of Claim 1 with
$r=3$ and $s=381$, we conclude that $m$ must have more than 381
distinct prime factors.

\def\ppp#1#2#3{p_{38#1}p_{38#2}p_{38#3}} 

Another computation reveals that the only numbers with at least 382
distinct prime factors that are less than $z$ are the numbers
$p_4p_5\dotsm p_{382}m'$, where $m'\in\{\ppp345$, $\ppp346$,
$\ppp356$, $\ppp347$, $\ppp357$, $\ppp456$, $\ppp348$, $\ppp367\}$;
and none of these numbers are congruent to $1\mod{30}$.
\vskip12pt

Let us define $n=(z-1)/30$, which by Claim 2 is both an integer and
the smallest possible solution of $\phi(30n+1)<\phi(30n)$. (Small
wonder that we haven't stumbled across any solutions of this
inequality---$n$ has 1,116 digits!) It would be quite gracious of $n$
to be an actual solution, and indeed it is.

First we show that $\phi(30n+1)/(30n+1)<\phi(30n)/30n$. We have
already computed that
\begin{equation}
{\phi(30n+1)\over30n+1} = \frac{\phi(z)}z = 0.2666117\dots.
\label{small}
\end{equation}
It turns out that $n$ is divisible by both 60 and
$p_{4,874}={}$47,279, so let us define $n'=n/(60p_{4,874})$. We can
compute that $n'$ is not divisible by any of the first 80,000
primes. This computation can be done most quickly by multiplying the
primes together in blocks of 1,000, say, and computing the greatest
common divisor of $n'$ and the product. Since computing greatest
common divisors is a very fast operation, checking that $n'$ is not
divisible by any of the first 80,000 primes takes only a few minutes
on a workstation---much more reasonable than trying to factor a number
with over a thousand digits.

Now define $q_2=\prod_{i=80,001}^{80,186}p_i$. We compute that $q_2$
has 1,118 digits and so $q_2>n>n'$. By using parts (a) and (b) of
Claim 1 with $r={}$80,000 and $s=186$, we see that
$\phi(n')/n'\ge\phi(q_2)/q_2$.  Therefore, since
$\phi(ab)=\phi(a)\phi(b)$ when $a$ and $b$ are relatively prime, we
compute that
\begin{equation}
{\phi(30n)\over30n} = {\phi(30\cdot60p_{4,874})\over
30\cdot60p_{4,874}} {\phi(n')\over n'} \ge \frac4{15} \big(
1-\frac1{47,279} \big) \frac{\phi(q_2)}{q_2} = 0.2666124\dots.
\label{big}
\end{equation}

This shows that $\phi(30n+1)/(30n+1)<\phi(30n)/30n$, which doesn't
quite imply that $\phi(30n+1)<\phi(30n)$, but only
$\phi(30n+1)<\phi(30n)(1+1/(30n))$. However, the numbers computed in
(\ref{small}) and (\ref{big}) differ in the sixth decimal place, while
multiplying by $1+1/(30n)$ leaves a number unchanged until past the
1100th decimal place.

Therefore we have proved:
\vskip12pt

{\noindent\bf Theorem. \it The smallest solution of
$\phi(30n+1)<\phi(30n)$ is}
\begin{equation*}
\openup-1\jot 
\begin{split}
\scriptstyle n=&\scriptstyle232,909,810,175,496,793,814,049,684,205,233,780,004,859,885,966,051,235,363,345,311,075,888,344,528,723,154,527,984,\\
&\scriptstyle260,176,895,854,182,634,802,907,109,271,610,432,287,652,976,907,467,574,362,400,134,090,318,355,962,121,476,785,712,\\
&\scriptstyle891,544,538,210,966,704,036,990,885,292,446,155,135,679,717,565,808,063,766,383,846,220,120,606,143,826,509,433,540,\\
&\scriptstyle250,085,111,624,970,464,541,380,934,486,375,688,208,918,750,640,674,629,942,465,499,369,036,578,640,331,759,035,979,\\
&\scriptstyle369,302,685,371,156,272,245,466,396,227,865,621,951,101,808,240,692,259,960,203,091,330,589,296,656,888,011,791,011,\\
&\scriptstyle416,062,631,565,320,593,772,287,118,913,728,608,997,901,791,216,356,108,665,476,306,080,740,121,528,236,888,680,120,\\
&\scriptstyle152,479,138,327,451,088,404,280,929,048,314,912,122,784,879,758,304,016,832,436,751,532,255,185,640,249,324,065,492,\\
&\scriptstyle491,511,072,521,585,980,547,438,748,689,307,159,363,481,233,965,802,331,725,033,663,862,618,957,168,974,043,547,448,\\
&\scriptstyle879,663,217,971,081,445,619,618,789,985,472,074,303,100,303,636,078,827,273,695,551,162,089,725,435,110,246,701,964,\\
&\scriptstyle021,045,849,081,811,604,427,331,227,553,783,590,821,510,091,607,567,178,842,569,576,699,548,038,217,673,171,895,383,\\
&\scriptstyle249,326,800,667,432,993,531,186,437,659,910,632,865,419,892,370,957,722,154,266,351,039,808,548,150,828,868,968,820,\\
&\scriptstyle675,198,820,381,135,523,646,361,202,383,915,218,571,017,801,463,011,491,108,784,343,253,284,393,511,650,254,506,597,\\
&\scriptstyle923,969,653,616,813,897,710,621,756,693,827,471,154,701,151,222,320,443,347,408,180,047,964,860.
\end{split}
\end{equation*}

\vskip12pt
I would like to thank Mike Bennett for verifying the above
computations and to acknowledge the support of National
Science Foundation grant DMS 9304580.

\vskip12pt\baselineskip=12pt\parindent=0pt\font\ten=cmti10\ten

School of Mathematics, Institute for Advanced Study, Olden Lane,
Princeton, NJ 08540, U.S.A.

gerg@@math.ias.edu

\end{document}